\documentclass{article}

\usepackage{enumerate}   
\title{A classification of one-dimensional local domains based on the invariant $(c-\delta)r-\delta$.} 
 
    \author{Anna Oneto$^{a }$ \ and \  Elsa Zatini$^b  ${\footnote{Corresponding author.
\newline \indent {\em E-mail addresses}:\quad  oneto@diptem.unige.it (A.Oneto), \quad  zatini@dima.unige.it (E.Zatini).}}
 \\
{\small $^{a }$ {\it        Ditpem, 
  Universit\`a di Genova,   P.le Kennedy, Pad. D -  I 16129 Genova (Italy);}} \\ {\small {\it$^b$
    Dima,  Universit\`a di Genova,  Via Dodecaneso 35 -  16146 Genova (Italy) }}}
\date{}

\expandafter\chardef\csname pre amssym.def  at\endcsname=\the\catcode`\@

\catcode`\@=11

\def\undefine#1{\let#1\undefined}
\def\newsymbol#1#2#3#4#5{\let\next@\relax
    \ifnum#2=\@ne\let\next@\msafam@\else
    \ifnum#2=\tw@\let\next@\msbfam@\fi\fi
    \mathchardef#1="#3\next@#4#5}
\def\mathhexbox@#1#2#3{\relax
    \ifmmode\mathpalette{}{\m@th\mathchar"#1#2#3}%
    \else\leavevmode\hbox{$\m@th\mathchar"#1#2#3$}\fi}
\def\hexnumber@#1{\ifcase#1 0\or 1\or 2\or 3\or 4\or 5\or 6\or 7\or 8\or
    9\or A\or B\or C\or D\or E\or F\fi}

\font\tenmsa=msam10
\font\sevenmsa=msam7
\font\fivemsa=msam5
\newfam\msafam
\textfont\msafam=\tenmsa
\scriptfont\msafam=\sevenmsa
\scriptscriptfont\msafam=\fivemsa
\edef\msafam@{\hexnumber@\msafam}

\font\tenmsb=msbm10 
\font\sevenmsb=msbm7
\font\fivemsb=msbm5
\newfam\msbfam
\textfont\msbfam=\tenmsb
\scriptfont\msbfam=\sevenmsb
\scriptscriptfont\msbfam=\fivemsb
\edef\msbfam@{\hexnumber@\msbfam}

      \font\tengothic=eufm10
      \font\sevengothic=eufm7
      \newfam\gothicfam
      \textfont\gothicfam=\tengothic
      \scriptfont\gothicfam=\sevengothic
      \def\goth#1{{\fam\gothicfam #1}}
      \font\tenmsb=msbm10
      \font\sevenmsb=msbm7
      \newfam\msbfam
      \textfont\msbfam=\tenmsb
      \scriptfont\msbfam=\sevenmsb
      \def\Bbb#1{{\fam\msbfam #1}}
\newtheorem{prop}{Proposition}[section]
\newtheorem{rem}[prop]{Remark}

\newtheorem{thm}[prop]{Theorem}
\newtheorem{coro}[prop]{Corollary}

\newtheorem{lemma}[prop]{Lemma}
\newtheorem{notat}[prop]{Notation}
\newtheorem{example}[prop]{Example}
\newtheorem{set}[prop]{Setting}
\newtheorem{(*)}[prop]{}

\newcommand{\integ}{\mbox{${\sf Z\hspace{-1.3 mm}Z}$}}

\newcommand{\nat}{\mbox{${\rm I\hspace{-.6 mm}N}$}}

\newcommand{\R}{\mbox{$\overline{R}$}}

\newcommand{\m}{\mbox{${\goth m}$}}

\newcommand{\w}{\mbox{$\omega$}}

\newcommand{\I}{\mbox{$\  \Longrightarrow \ $}}
\newcommand{\II}{\mbox{$\  \Longleftrightarrow \ $}}

\newcommand{\enu} {\begin{enumerate}[{\rm(1)}]} 
\newcommand{\enua} {\begin{enumerate}[ $(a)$]} 
\newcommand{\denu} {\end{enumerate}}

\font\tengothic=eufm10
\font\sevengothic=eufm7
\newfam\gothicfam
\textfont\gothicfam=\tengothic
\scriptfont\gothicfam=\sevengothic
\def\goth#1{{\fam\gothicfam #1}}

\begin{document}

  \maketitle

 {\bf Abstract}.{ \small \  Let $(R, \m)$ be a one-dimensional, local, Noetherian domain
and let $\R$ be the integral closure of $R$ in its quotient field
$K.$ We assume that $R$ is not regular, analitycally irreducible and residually
rational.   The usual valuation $v : K\longrightarrow \integ \cup \infty $
associated to $\R$ defines  the numerical semigroup
  $v (R)=\; \{ v (a),\;
a\in R,\; a\neq 0\} \subseteq \nat $.
  The aim of the paper is to study the non-negative invariant
 $b:=(c-\delta)r- \delta $, where $c, \delta, r$ denote the conductor,  the length of $\R/R$ and  the  Cohen Macaulay type  of $R$, respectively. In particular, the classification of  the   semigroups $v(R)$ for rings  having $b\leq 2(r-1)$ is realized. This method of classification might be successfully utilized with similar arguments  but more boring computations in the cases $b\leq q(r-1), $ for reasonably low values of $q$.  The main tools are   type sequences and the invariant  $k$ which estimates the number of elements in $v(R)$ belonging to the interval $[c-e,c), \ e$ being the multiplicity of $R$. }\\

{\bf  Introduction. }
 Let $(R, \m)$ be a one-dimensional, local, Noetherian domain
and let $\R$ be the integral closure of $R$ in its quotient field
$K.$ We assume that $R$ is not regular and analitycally irreducible, i.e. $\R$ is a DVR with  uniformizing parameter $t$ and a finite $R$-module.
 We also suppose $R$ to be residually
rational, i.e. $R/\m
\simeq \R/t\R.$  
Called $v : K\longrightarrow \integ \cup \infty $   the usual valuation
associated to $\R $, 
the image $v (R)=\; \{ v (a),\;
a\in R,\; a\neq 0\} \subseteq \nat $ is a numerical semigroup. 
Starting from  the following classical invariants: \par
    $c  $,\   the {\em conductor} of $R,$  i.e. the minimal
$j\in v (R)$ such that   $  j+\nat \subset  v (R),$\par   
   $\delta :=  \ell_R (\R/R) $, the number of gaps of the semigroup $v (R)$ in
$ \nat ,$
\par $r:= \ell_R((R:\m)/R)$,   the {\em Cohen Macaulay type} of $R$,\\ 
the new invariant\\
\centerline{$b:=(c-\delta)r- \delta  $}   has been recently considered in the literature. The general problem of classifying rings according to the size of $b$  has been   examined by several authors.  First,
Brown and Herzog in   
\cite{bh} characterize   all the one-dimensional reduced local rings having $b=0$ or $b=1$.    Successively,   in  \cite{dd}, 
        \cite{dm},   \cite{dlm}, Delfino, D'Anna and Micale   consider  the   rings for which  $b\leq r$.  Partial answers in the
case 
$b>r-1$ are given in  [5]. 
 
  In   \cite[Section 4]{ozinv} we obtain some improvements  of the quoted results. This is done  by using the expression of the invariant $b$ in terms of the
type sequence
$[r_1,..,r_n]$ \  (defined in (\ref{set2})),   where $n:=c-\delta$ and $r_1$ equals  the Cohen-Macaulay type $ r $ of $R$,   namely: \\
  \centerline{$ b=\sum_{i=1}^n (r-r_i). $}
 So, employing the properties of the type sequence, we get as  a straightforward   consequence of the preceding formula the well known  bounds \\
 \centerline{$0 \leq b\leq (n-1)(r-1)$ }
 (for the positivity see \cite{bh}, Theorem 1;  for the upper bound see \cite{dd}, Proposition 2.1). Also, we recover in an immediate way the two extremal cases:\par
 $b=0$, corresponding to the so called rings {\it of \ maximal \ length}, i.e. the rings having   maximal type sequence $[r,r,...,r]$; \par $ b= (n-1)(r-1)$, corresponding to the   {\it almost \ Gorenstein} rings, i.e. the rings having   minimal type sequence $[r,1,...,1]$.\\
Actually,  for any integer $q\in\nat$ it is natural to ask if it is possible to characterize  the rings     verifying \\
\centerline{$(q-1)(r-1)   \leq b \leq q(r-1).$} In Section 3   we write explicitly all the possible values of $v(R)$  for   $1\leq q\leq 2$  (see Theorems  (\ref{b=r-1}), (\ref{b-compreso}), (\ref{b=2(r-1)})), but we  outline that the method used here is absolutely general and analogous although more tedious computations might be repeated for greater values of $q$.    \\         
To achieve our results, we utilize heavily  the number  
  \\ \centerline{$k:=\ell_R (R/({\goth C}+xR)),$} 
 where ${\goth C}:=t^c \R$ denotes the $ conductor \ ideal$ of $R$ in $\R$ and $x$ an element of $ R$   such that $v(x)=e(R)$, the {\it multiplicity}. In   \cite{d}   it is established that $b=r-1 \I   k=1$ or $2$  \cite[Proposition 2.4]{d}, and that $ b=r-1 $ and 
$k = 2 \I r=e-2$  \cite[Corollary 2.13]{d}.  In \cite{dlm}  
the lower bound $r k -e+1\leq b $ is found. Improvements  of these results and  several other  inequalities relating the  invariants $k, b, r$ are now realized by means of the
  type sequence of $R$ (see (\ref{dopoD1}) and   (\ref{dopoD})).  For this purpose we   introduce in Section 1  a decomposition of the semigroup $v(R)$  as 
a disjoint union of   subsets:
\\ \centerline{$v(R)\!=\!   \{0, e,2e,...,pe, c,\rightarrow\} \cup H_1 \cup .... \cup H_{k-1},$}  where  \
 $ H_i:=\{y_{i},y_{i}+e,...,y_{i}+l_{i}e\}$,  \  $i=1,...,k-1, \ \ p, \ l_i \in \Bbb N,$ and  $\{y_i\}_{i=1,..,k-1}$ have distinct residues (mod $e$)  (see Setting
\ref{set3}).  This
allows us to obtain in Section 2
 the   useful  formula (\ref{XYZ}.1):\\
\centerline{$ b= X+ Y +Z$}  
 where  $\ X:=(k-1)(r-1)\geq 0$, \par
 $ \ \ \ \  \ \! \   Y:=k-(e-r) \geq 0$, \par $   \ \ \ \ \ Z:=(r+1)(p+\sum_1^{k-1}  l_i)+k+h-pe-1\geq 0$.\\
 Obviously $X+Y=rk-e+1$, and so the integer $Z$ measures how far is $b$ from the   lower bound  proved in \cite{dlm}.   \\ 
 The advantage of this formula is evident for low values of $b$. For instance,  for rings having $b\in\{0,1,2\}$
  we state in a  quite simple way all  the possible value
sets  (see Theorems (\ref{b=0}), (\ref{b=1}), (\ref{b=2})). Nevertheless,
 a such type of classification  might  be  accomplished for greater values of $b$ with similar arguments.  
 
\section{Preliminary results.}
\hspace{0.3cm} We begin by giving the setting of the paper.

\begin{set} \label{set2}  {\rm    Let ($R, \m$) be a one-dimensional local Noetherian domain with  residue field  $\kappa$ and quotient field $K$. We assume
throughout that  $R$ is not regular with  normalization
 $\overline{R} \subset K$   a DVR and a finite $R$-module, i.e., $R$ is analytically irreducible.
Let $t\in\R$ be a uniformizing parameter for $\R$, so
that $t\R$ is the maximal ideal of $\R$. We also suppose that the field $\kappa$ is  isomorphic to the residue field   $\R/t\R$, i.e., $R$ is residually rational. We denote the  usual   valuation on $K$ 
associated to $\R $
by $v$; that is,   $ v : K\longrightarrow \integ \cup \infty$,  and $v(t)=1$. 
By \cite[Proposition 1]{mat} in this setting it is possible to compute   for a  pair of fractional nonzero ideals   $  I
\supseteq J $  the length of the $R$-module $I/J$  
  by means of valuations:  \\ \centerline{ 
  (\ref{set2}.1) \hspace{1.5cm}
  $\ell_R (I/J)= | v(I) \setminus v(J)|. $\hspace{3.5cm}}  
       The set  $v(R)  :=   \{v(a)\ | \ a\in R,\; a\neq 0\} \subseteq \nat $   is the {\it   numerical semigroup} of
$R$. Since the {\it conductor}  ${\goth C}:=(R:_K \!\R)$ is an ideal of both $R$ and $\R$, there exists a positive integer $c$ so that $ {\goth C}   =t^c\overline R, \ \ell_R(\R/{\goth C})=c$ and $c \in v(R)$. Furthermore, denoting by    $\delta :=  \ell_R (\R/R) $  the number of gaps of the semigroup $v (R)$ and  $r:= \ell_R((R:\m)/R)$    the {\em Cohen Macaulay
type} of
$R$, we define the invariant \\
\centerline{$b:=(c-\delta)r- \delta $.} We list   the elements of   $v(R)$ in order of
size:
$v(R):=\{s_i\}_{i
\geq 0}$, where  
 $s_0=0$ and $s_i< s_{i+1} $, for every $i\geq 0$. We put $e:= s_1$ the {\it multiplicity} of $R$ and $n=c-\delta$ the number    such that $s_n=c$. For every $i \geq 0$, let $R_i$ denote the ideal of elements whose values
are bounded by $s_i$, that is,  \\
\centerline{$R_i :=  \{ a\in R\; |\; v(a)\geq s_i\}.$}
The  ideals  $R_i$    give a strictly decreasing sequence\par $R=R_0\supset R_1=\m\supset R_2\supset
\ldots \supset R_n ={\goth C} \supset R_{n+1}\supset...\ ,\vspace{0.1cm} $\\
which induces  the chain of duals:  \par 
  $R\; \subset \; (R\; :R_1)\; \subset \; 
...  \subset (R\; :R_n)=\; \R \; \subset (R\;
:R_{n+1})=t^{-1}\R\subset...    $.\\
Put  \ $r_i:=l_R((R:R_{i})/(R:R_{i-1})), \  i\geq 1;$   the finite sequence of  integers \centerline{$[r_1,\dots,r_n] $ \quad is  
the {\it type sequence}     of   $R$.} 
In particular $r_1=r,$    the  Cohen-Macaulay type    of $R$.  Moreover it is known that:\par
 $\bullet\quad 1\leq r_i \leq  r  $  for every $ i\geq 1$, \ and \   $  r_i=1 $ for every $  i>n$, \par     $ \bullet\quad  \delta   = \sum_1^n \  r_i, $ \par
    $\bullet\quad2\delta -c =  \sum_1^n  (r_i-1) = \sum_1^{\infty} (r_i-1)$ (see, e.g. \cite[ Prop.2.7]{ozinv}).}

\end{set}

 \noindent  Type sequence is a suitable tool   to study the behavior of the invariant  $b$.   \begin{prop}\label{n-1}  We have:
\enu \item $ \    b=\sum_{i=1}^n (r-r_i).   $ \ \item \  $0\leq
b\leq(n-1)(r-1).\qquad$ 
\denu
 \end{prop}  Proof.  For (1) see \cite[Section 4]{ozinv}. \\
 (2).   We have:  $   \sum_{i=1}^n (r-r_i)  = \sum_{i=2}^n (r-r_i )\leq (n-1)(r-1)$, because
$r_1=r$ and 
$r_i\geq 1$, \ for every $i\geq 1$.\quad $\diamond$

\begin{notat}\label{notaB}  {\rm   Let $R$ be as in  (\ref{set2}). We set:  \enu
 \item[$\bullet$] $ x \in \m $ is  an element   such that
$v(x) =e$;  namely, $\ell_R (R/ xR)=e$.  
 \item[$\bullet$]  For $a,b \in \Bbb Z$, \ \  $[a,b]= \{ x \in  \Bbb Z \ | \ a \leq x \leq b\}$. 
\item[$\bullet$]    $i_0 \in [1,n] $ is such that $s_{i_0-1 } = $ min$\{ y\in v(R)  \ |  \ y\geq c-e\}.$\\
$(i_0=1 \iff c=e)$.
\item[$\bullet$]   $B:= [i_0,n]$ \  and \    $A:= [1,n] \setminus B$ \ \ ($|A|\leq n-1)$. 
\item[$\bullet$]  $k:=\ell_R (R/({\goth C}+xR)) $ \ \ $(1\leq k\leq e-1)$.
   \end{enumerate}}\end{notat}
 
\begin{thm}\label{dopoD1}  The following facts hold.
\enu
 \item $k=|B|= \ell_R ({\goth C}:_R \m/ {\goth C})  \geq e-r>0$.
 \item  $ k\leq\sum_{i\in B} r_i \leq e-1. \ \ If $   
   $ \sum_{i\in B} r_i= e-1,$  then  $  s_{i_0-1}=c-e$. 
\end{enumerate} \end{thm}

\underline{Proof}. (1) and the inequality 
$\sum_{i\in B} r_i
\leq e-1
$ of (2) are proved in \cite[Lemma 4.2]{ozinv}. Since $r_i \geq 1 $ for every $i$ and $|B|=k$, the inequality $ k\leq\sum_{i\in B} r_i $ is done.  \\
  Moreover, denoting by $\w$ the canonical module of $R$ (see \cite{ozinv} for the existence and the properties in our setting),  we  remark  that\par 
  $ \sum_{i\in B} r_i= \ell_R (\R/(R:R_{i_0-1}))= |v(\w R_{i_0-1})_{<c} |$ and \par $v(\w R_{i_0-1})_{<c}  \subseteq [ c-e, c-2]$  \\
  (see the proof of the quoted lemma).
 Thus    $ \sum_{i\in B} r_i= e-1 \I v(\w R_{i_0-1})_{<c} = [ c-e, c-2]$,\  and so $  s_{i_0-1}$, the minimal element in $ v(\w R_{i_0-1})$,
equals $c-e$.
\quad $\diamond$ \\

 The case $k=1$ is
completely known and recalled below  for the convenience of the reader.   
\begin{prop}\label{delta1}    {\rm \cite[Lemma 4.4]{ozinv}}
  The following facts are equivalent:
\enu
\item    $k =1.$
\item   $v(R)=\{0,e,...,pe,c\rightarrow\}.$
\item The type sequence of $R$   equals $[e-1,....,e-1,r_n]$.
\end{enumerate} 
If $R$ satisfies these equivalent conditions, then:   \enu \item[]
$  \delta=c-p-1, \   b=(p+1)e-c\leq r-1,\ r=e-1, \  r_n=e-1-b. $
\denu
\end{prop}

 By virtue of    (\ref{set2}.1) we have   $ k= |v(R) \setminus v({\goth C}+xR)|$. This fact allows   to  write  $v(R) = v( {\goth C}+xR )
\cup
\{0, y_1, ..., y_{k-1}\} $, obtaining   the description of $v(R)$ as a disjoint union of the sets $H_i$ given in the next setting. 
The construction is significant for $k >1$.  

\begin{set} \label{set3}  {\rm    Let $k  >1$. We set: \enu \item[] $v(R)\!=   \{0, e,2e,...,pe, c,\rightarrow\} \cup H_1 \cup .... \cup H_{k-1}$,  
where 
  \item[$\bullet$]    $p  $ is  the
integer such that \  $c-e\leq pe<c $, in other words, 
$pe+2\leq c\leq (p+1)e$.  \ \ $(p\geq 0$ and $\ p=0 \iff c=e)$.
 \item[$\bullet$]   $h:= (p+1)e-c $, \ \ \ $ (0\leq h\leq e-2).$
 \item[$\bullet$] $ H_i:=\{y_{i},y_{i}+e,...,y_{i}+l_{i}e\}$,  \  $i=1,...,k-1, \ l_i \in \nat.$ 
\item[$\bullet$] The integers  $ y_i  \in \nat$ are   such that 
$  e < y_1<y_2< ...< y_{k-1}$, \  $y_i \notin e \Bbb Z$, 
 \    $ \overline{y_i} \neq \overline{y_j} \  (mod \ e)  $ for every $\ i,j\in \{1,..,k-1\}$.
 \item[$\bullet$]The integers $l_i , \    i=1,...,k-1$, are defined by the relations:\par
$ y_{i}+l_{i}e <c \leq y_{i}+(l_{i}+1) e$.
 \item[$\bullet$] For $ k=  2 $ we shortly call \ $y:=y_1, \ \ l:=l_1$.
\denu}
     \end{set}
     
     \begin{example} {\rm If $v(R)=<10,11,26>$, then:\\
    $v(R)=\{0,10,20,30,40,50 \to\} \cup H_1\cup ... \cup H_7$ where $H_1=\{11,21,31,41\}, \ H_2=\{22,32,42\}, \ H_3=\{26,36,46\}, \ H_4=\{33,43\}, \ H_5=\{37,47\}, \ H_6=\{44\}, \ H_7=\{48\}$.  According to notations previously introduced $y_1=11, y_2=22, y_3=26,y_4=33,y_5=37,y_6=44,y_7=48$ and $l_1=3, l_2=l_3=2, l_4=l_5=1, l_6=l_7=0$. Moreover, $c=50, p=4, h=0$.   }\end{example}
     
\begin{prop} \label{aftersetting}  Let    $k>1, p, h, \{l_i\}  $   be the
integers defined in $(\ref{notaB})$ and $(\ref{set3})$. Then:
				\enu
\item $r \in \{ e-k, ...., e-1\}.$
\item    $ 0\leq l_{k-1} \leq .... \leq l_2 \leq l_1\leq p-1 $.
 \item $ c-\delta= p+k+\sum_1^{k-1} l_i$,\par
 $ \delta= (p+1)(e-1)-h -\sum_1^{k-1} (l_i+1)$.
 \end{enumerate}
 \end{prop} 

\underline{Proof}. Assertion (1)  follows
immediately from    (\ref{dopoD1}.1). \\
 (2). By definition of $l_i$ and $p$,  we have   \
 $(l_i+1)e< y_i+l_ie <c \leq (p+1)e $; \\ then
  $l_i+1\leq p$, \ for every $i=1,...,k-1.$ \ Now  note that \par
 $  y_i+l_ie <c \leq y_{i-1}+(l_{i-1}+1)e \I y_i-y_{i-1}<(l_{i-1}+1-l_i)e \I l_i \leq l_{i-1}$.\\            
(3). Using the integers defined in (\ref{set3})       $c- \delta$ and $\delta$ can be expressed as :\par
 $c- \delta=(p+1)+(l_1+1)+...+(l_{k-1}+1)=p+k+\sum_1^{k-1} l_i$,\par
 $\ \ \ \  \ \delta =c-(c- \delta)= (p+1)e-h -(p+k+\sum_1^{k-1} l_i)$\par
 \qquad\   $=(p+1)(e-1)-h -\sum_1^{k-1} (l_i+1)$.  
\quad 
$\diamond$ \\

It is natural to ask how the   elements  $y_1,...,y_{k-1}$ introduced in (\ref{set3})  influence the Cohen Macaulay type of $R$.    This will be analysed in the following
(\ref{tipoe-1}), (\ref{l=k}), (\ref{l=2}).

\begin{prop} \label{tipoe-1} Let $k=\ell_R (R/({\goth C}+xR)) $ and let $v(R)$ be as in $(\ref{set3})$.   Further let
$x_1,...,x_{k-1} \in \m$ be   such that $v(x_i)=y_i$. The following facts are equivalent: 
				    	\enu
\item $ r=e-1 $, i.e. $R$ is of maximal Cohen Macaulay type.
 \item $ v(R) \setminus v(xR:\m) =\{0 \}.$  
\item $y_1,...,y_{k-1} \in v (xR: \m)$.
\item $x_1,...,x_{k-1} \in  (xR: \m)$.
\item $x_ix_j \in x\m $ for every $i, j = 1,...,k-1$.
\item $\ell_R(\m /\m^2)=e$,  i.e. $R$ is of maximal embedding dimension.
 \denu
   \end{prop}
 
\underline{Proof}. Since
  $e-r=\ell_R(R/xR)-\ell_R( (xR:\m) /xR)= \ell_R(R/  (xR:\m))$, the equality $e-r=1$ means $|v(R) \setminus v(xR:\m)|=1$, and so $1 \iff 2$ is proved. In the same way we obtain that \par
   $(*) \ \ \ r=e-1 \iff (xR:\m) = \m \iff   \m^2 = x \m $.  \ Moreover, \par
    $(**) \ \ v(x^{-1}\m ) \subseteq \Bbb N \I  x^{-1}\m  {\goth C}\subseteq {\goth C} \I \m {\goth C} = x {\goth C}\I  {\goth C}\subseteq (xR: \m)$.\\
  Considering the chain of ideals \par
  $R \supset \m \supseteq {\goth C} +(x,x_1,...,x_{k-1})R \supset  {\goth C} + (x,x_1,...,x_{k-2})R  \supset ...  \supset {\goth C}+xR$, 
 \\we see that   $ \ell_R (R/({\goth C}+xR)) =k\I \m= {\goth C} +(x,x_1,...,x_{k-1})R $, hence \par
 $(***) \ \ x_i\m= (xx_i)R+(x_ix_j)R+x_i{\goth C}$ for  every $  j=1,...,k-1 $.\\
By  $(*) $   we have immediately $1 \iff 6$ and $1 \I 5$.\\
$5 \I 4$. By the assumption  $x_ix_j \in x\m, \ \forall \ i, j = 1,...,k-1$ and by the obvious inclusion $x_i{\goth C} \subseteq \m  {\goth C}=x{\goth C}$, from
$(***)$  we get $x_i\m  \subseteq xR$, then $x_i \in (xR: \m)$. 
 \\
The implication  $4 \I 3$ is   obvious.\\ Finally,   $3 \I 2$ holds by $(**)$.
\quad $\diamond$

\begin{rem} \label{es1} {\rm It is clear  from (\ref{tipoe-1}) (see the equivalence 1-5) that   the condition  
  $y_i+y_j-e \in v(R) \ {\rm for \ every } \ i, j=1,...,k-1,$ 
is necessary to have maximal Cohen Macaulay type. Unfortunately,  it is not sufficient. For example,
 if $R=\kappa[[t^{6},t^{9}+t^{10}, t^{14},t^{16},t^{17},t^{19} ]]$, then $k=2, \ y=9, \ 2y-e=12  \in v(R) $, but $ \ r=e-2 $. In this case (\ref{tipoe-1}.5) does not hold, because $(t^{9}+t^{10})^2\notin xR$.}
\end{rem}

\begin{prop} \label{l=k}  Let $k=\ell_R (R/({\goth C}+xR)) $ and let $v(R)$ be as in $(\ref{set3})$. 
				\enu
				   \item   
       $ r=e-k, \ k \geq 2, \II  v(R) \setminus v(xR:\m) =\{0,y_1,....,y_{k-1}\}.$ 
 \item If $ r< e-1,$  then \vspace{-0.3cm}
       \enua
 \item $2y_1   <c+e,$  
 \item $p \leq  2l_1+ 2$  and $p= 2l_1+ 2\I h>0.$
       \denu
    \item If $ r=e-k,$   then \vspace{-0.3cm}
             \enua
 \item $y_1+y_j  <c+e,$   for every  $  j=1,...,k-1.$
 \item $p \leq  l_1+l_{k-1}+2  $   and $ p= l_1+l_{k-1}+2 \I h>0.$
           \denu
   \item   If $  p\geq 3 $ and $ i$ is such that $l_i=0$, then $  2y_i > c+e.$
\end{enumerate}
 \end{prop}
 
\underline{Proof}. (1).  By means of $(**) $ stated in the   proof of (\ref{tipoe-1}), we have the inclusions  
  $({\goth C}+xR) \subseteq (xR:\m) \subseteq R$. Since   $e-r= \ell_R(R/  (xR:\m))$ and
$k= \ell_R(R/({\goth C}+xR))$, it follows   that  
$e-r=k \iff ({\goth C}+xR) = (xR:\m)$.\\
 To see  (2.a), suppose   
$2y_1   \geq c+e,$  then   $y_i+y_j   \geq c+e$   for every  $i,j=1,...,k-1.$
   Let
$x_{i}\in \m$ be    elements such that $v(x_{i})=y_{i}$ and let $s \in \m$. If $s \in ({\goth C}+xR)$, then $x_i s \in \m ({\goth C}+xR)\subseteq xR$.  If $s \notin ({\goth C}+xR)$, then $v(s)=y_j,$ for some $j,  \ 1\leq j \leq k-1$, hence $v(x_i s)=y_i+y_j \geq c+e \I x_i s \in x{\goth C} \subset
xR$. In both cases $x_i \in (xR:\m)$, and so $y_i \in v(xR:\m)$. Thus   $v(R)
\setminus v(xR:\m) =\{0\}$ and $r=e-1$ by (\ref{tipoe-1}), a contradiction.\\ To see (2.b),  consider that by (\ref{notaB}) and (\ref{set3}): \par   
   $ y_1 \geq c-(l_1+1)e= (p-l_1)e -h$. \\
   Combining this with the preceding   (2.a),      we obtain \par
 $(2p-2l_1 )e -2h\leq 2y_1 < c+e=(p+2)e-h $. \\  Thus \    $ (p-2l_1 -2)e< h$ and since   $h \leq e-2$,    we  see that $ p\leq  2l_1 +2 $ and also  
 that   $\ p=2l_1+ 2 \I h>0.$ \\
To prove (3.a), it suffices to show that $y_1+y_{k-1}  <c+e.$ Suppose    $y_1+y_{k-1} \geq c+e $, then  $y_i+y_{k-1} \geq c+e $ for all $i$. Let $x_{k-1} \in \m$ be an element such that $v(x_{k-1})=y_{k-1}$.
As in  (2.a), we get $x_{k-1} \in (xR:\m)$, and so  $ y_{k-1}\in v(xR:\m)$,  a contradiction, since  the assumption   $ e-r=k $   means  $  v(R)
\setminus
v(xR:\m) =\{0,y_1,....,y_{k-1}\} $ (see item 1). \\
 We prove now (3b). As in  (2.b), \par $ y_j\geq c-(l_j+1)e= (p-l_j)e -h$, \  for $j=  1,..., k-1,$ \\ and by (3.a) \par  
$(2p-l_j-l_1)e-2h\leq  y_1+y_j<c+e= (p+2)e-h$. Hence \par $(p-l_j-l_1-2)e<h \leq e-2$, for every $j=1,...,k-1$.\\
We conclude \par $p\leq l_1+l_j+2 \leq l_1+l_{k-1}+2 $ and also the last assertion.\\
For (4), note that $l_i=0  \I y_i+e \geq c$, and that $p\geq 3  \I c> 3e$.\\ Thus:  \  $2y_i \geq 2c-2e =c+ (c-2e)> c+e $, as desired.  \quad $\diamond$  \\

We may describe the particular  case $k=2$ in a more precise way.  

 \begin{prop}\label{l=2}  Assume  $ k=  2.$   With setting  {\rm (\ref{set3})} we have: 
  \enu
  \item  $r=e-1 \II$ one of the following conditions is satisfied:
\enua \item  $   \ 2y \geq c+e;$ \item \
 $ 2y=(2q+1)e<c+e,   \ q\geq 1,  \   p \geq 2$  and $\ y\in v(xR:\m)$.
\denu
    \item $r=e-2 \II 2y < c+e$ and if  \   $   2y=(2q+1)e $, \  then $\  y\notin v(xR:\m)$.
  \end{enumerate} 
 \end{prop}

\underline{Proof}.  
First recall that by (\ref{dopoD1}.1) one has $r\geq e-2$. For  implication $\I$ in (1), note that  $y\in v(xR:\m)$, by (\ref{tipoe-1}),  and  so $ 2y-e \in v(\m)$. Then regarding the
structure of
$v(R)$, we have the claim. For the opposite implication, note that in case ($a$) for any $s\in \m$ such that $v(s)=y,$  $v(x^{-1}s^2 )=2y-e \geq c \I x^{-1}s^2  \in  {\goth C} \I s^2 \in x
\m$; now use again  (\ref{tipoe-1})  to conclude.   
 \\ (2) is immediate by (1).\quad $\diamond$

 \section{ Bounds for the invariant $b$.} 
 
Starting from the preliminary result (\ref{n-1}) we go on in studying the integer $b$.
First (see (\ref{dopoD})) we find lower and upper bounds using the properties of the type sequence, then  (see (\ref{XYZ})) we  express $b$ in terms of the integers $k, p, l_i, h$ occurring in the decomposition of $v(R)$  as in  (\ref{set3}). This description becomes quite simple in the particular cases $k=2,3$ (see (\ref{Casek=2}) and (\ref{Casek=3})). The last result of   the present section  (see(\ref{aggiunto1})) furnishes  informations according to the range  $\!(q-1)(r-1)< b\leq \! q(r-1)\!$, that will be basic in the next section.
 
\begin{prop}\label{dopoD}   With  Notation $\ref{notaB}$, the following facts hold.
\enu
\item   $  (e-r-1)(r-1)   \leq  rk-e+1\leq b-\sum_{i\in A}(r-r_i )  \leq k(r-1) $.   
 \item  
 $b=(k-1)(r-1)  + \sum_{i\in A}(r-r_i )\II   \sum_{i\in B} r_i= e-1$ and $k=e-r$. 
 \item $b=k(r-1)+\sum_{i\in A}(r-r_i )\II r_i=1  $ for every   $ i\in B$. 
\item The following conditions are equivalent: \enua \item  $b= (e-r-1)(r-1)$.\par
\item $b= (k-1)(r-1)$.
\item $  e-r=k,  \ \sum_{i\in B} r_i= e-1$   and $r_i=r  $ for every   $ i\in A $.\
\denu
 If these conditions  hold, then   $s_{i_0-1}=c-e$.
\item $b \geq (r-1)s$, where $s:=| \{ i \in [1,n] \ | \ r_i=1\}| $.
 \end{enumerate} \end{prop}

\underline{Proof}. 
 Write the invariant $  b=\sum_{i=1}^n (r-r_i )$ in the following form:\vspace{-0.2cm}$$(*)  \qquad b= \sum_{i\in B}(r-r_i )+ \sum_{i\in
A}(r-r_i )   
  =rk- \sum_{i\in B} r_i + \sum_{i\in A}(r-r_i ).\vspace{-0.2cm}$$  Using that $\sum_{i\in B} r_i\leq e-1$ (see
\ref{dopoD1}.2),   we obtain \\
\centerline{$ (**) \qquad rk- (e-1)  \leq b- \sum_{i\in A}(r-r_i )\leq k(r-1) $.\hspace{2cm}} \\ Then, since   $k\geq e-r$ by (\ref{dopoD1}.1),  the inequalities of (1) are clear.\\ 
(2). Supposing  \ $b-\sum_{i\in A}(r-r_i )=(k-1)(r-1)  $ we have by item 1 $(k-1)(r-1) \geq rk-e+1$, hence $k \leq e-r$ and since always $k\geq e-r$, it follows that $k=e-r$. From $(*)$  $ \sum_{i\in B} r_i= rk-(k-1)(r-1) =k+r-1=e-1$. For the converse, it suffices to substitute $ \sum_{i\in B} r_i= k+r-1$ in $(*)$.\\
(3). Using $(*)$ we have $b-\sum_{i\in A}(r-r_i )=k(r-1)  \II  \sum_{i\in B} (r-r_i)=k(r-1)$. Since $r-r_i \leq r-1$ for every $i$ and $k=|B|$, the last fact is equivalent to say that  $ r_i=1  $ for every   $ i\in B$.\\ 
(4), $a \I b$. By (1) we have  immediately $\sum_{i\in A}(r-r_i )=0$ and $  (e-r-1)(r-1)  =  rk-e+1 \I e-r=k$, as desired.\\ 
$b \I c$. By (1) we have $\sum_{i\in A}(r-r_i )\leq b-(rk-e+1)= -k-r+e \leq 0$, then we can apply item 2 with  $\sum_{i\in A}(r-r_i )=0$.\\
$c \I a$. Substitute in $(*)$ the relations of (c).\\The
fact   
$s_{i_0-1}=c-e$ is
  immediate by  (\ref{dopoD1}.2).   \\
  By applying \cite[Corollary 3.13.2]{ozinv}, with $I=\goth C$ \ we get (5). \quad $\diamond$    \\

Utilizing the description of the value set $v(R)$   introduced in (\ref{set3}),  we obtain the next useful formula for the invariant $b$.

\begin{thm} \label{XYZ} \ With Setting \ref{set3}, assume $k>1$. The following equalities hold: \enu 
\item  $ b=(r+1)  \sum_1^{k-1}  (l_i+1)- (p+1)(e-r-1)+h  
 =X+ Y +Z$ \par
 where  $X:=(k-1)(r-1)\geq 0$, \par \ \ \ \ \ \ \ \  $Y:=k-(e-r) \geq 0$,  \par \ \ \ \ \ \ \ \ $Z:=(r+1)(p+\sum_1^{k-1}  l_i)+k+h-pe-1\ \geq \sum_{i\in
A}(r-r_i )\geq 0$.
\item $c= (p+1+ \sum_1^{k-1}  (l_i+1))(r+1)-b$.

\denu \end{thm}

\underline{Proof}. (1). To get the desired formula it suffices to substitute in 
the equality $b=( c-\delta)r-\delta$ the expressions of $ c-\delta$ and $ \delta$ given in (\ref{aftersetting}.3). The
positivity of $ Y$ is clear by (\ref{dopoD1}.1). To prove the positivity of $Z$ we use the second inequality of (\ref{dopoD}.1):   $X+Y= kr-e+1\leq
b-\sum_{i\in A}(r-r_i )$, and so we have the conclusion:  $Z=b-(X+Y)\geq \sum_{i\in A}(r-r_i )\geq 0$. \\
(2). Since $b+c=(r+1)(c-\delta)$, (2) follows easily.\quad
$\diamond$  

 \begin{lemma}\label{Casek=2} {\rm Case} $k=2$. With Setting \ref{set3}, assume $k=2$. 
\enu
      \item    If $r=e-1$, then $\ b= (l+1)e+h \leq (l+2)e-2$.  \par
 Further: \ \ 
  $b =(l+2)e-2 \iff h=e-2$.
 \item If $r=e-2$, then:  \par
  \ \ \ $  b=(l+1)(e-1)+h -p-1 $,\par
  \ \ \  $c=(p+l+2)(e-1)-b. $\par Further we have:  \enua \item  
  $l+1\leq p \leq 2l+2 $  \  and \ $\ p=2l+2  \I h>0.$ 
 
 \item 
  $(l+1)(e-3)\leq b \leq(l+1)(e-2)+e-3$.  \ In particular 
\par
  $b=(l+1)(e-3) \iff p= 2l+2, h=1$ or $ p=2l+1, h=0$. \par  $ b =(l+1)(e-2)+e-3 \I p=l+1, \ p>1, \ h=e-2, \ y=e+1$.  
 \end{enumerate} 
\denu
 \end{lemma}
 
\underline{Proof}. For $k=2$,   we write      $v(R)\!=\!   \{0, e,2e,...,pe, c,\rightarrow\} \cup  
\{y ,y +e,...,y +l e\} $,  with $r \in \{ e-2, e-1\}, \ c-\delta= p+2+ l, \ c= (p+1)e-h$.    (see \ref{aftersetting}), (\ref{set3})). Then 
the expressions of $b, $ in items 1,2,  come from  (\ref{aftersetting}.4)  with   $k=2$ and $e-r=1, 2$, respectively. To complete the proof of item 1 recall that $h
\leq e-2$. \\
 The bounds for $p$ in item 2 come from (\ref{aftersetting}.2) and (\ref{l=k}.3) and the value of $c$ comes from (\ref{XYZ}.2). \\ Rewriting $b$ in the form \par
 $b=(l+1)(e-2)+(l-p)+h$,\\
 and recalling that $l-p \leq -1, \ h \leq e-2$, we obtain the upper bound for $b$.\\
 Rewriting $b$ in the form \par
 $b=(l+1)(e-3)+(2l+2-p)+(h-1)$, \\
 and using part $a$, we obtain the lower bound and also $  b=(l+1)(e-3)
\iff  p= 2l+2, h=1$ or
$ p=2l+1, h=0$.  \\  Finally, note that $ b =(l+1)(e-2)+e-3  \I h=(p - l-1)+e-2 \geq e-2 \I p=l+1, \ h=e-2 \I c=pe+2$ and since by definition of $l$  $ y+le<c$, it follows that $y<e+2$, hence $y=e+1$ and
$p>1$.   
\quad
$\diamond$

\begin{lemma}\label{Casek=3} {\rm Case} $k=3$. With Setting \ref{set3}, assume $k=3$. 
\enu
\item If $r=e-3 $, then $b=(l_1+l_2+2)(e-2)+h-2(p+1).$ Moreover,\par
 $p < l_1+l_2+2 \I b \geq (l_1+l_2+2)(e-4)+h,  \ h \geq 0.$ \par
  $p=l_1+l_2+2 \I b = (l_1+l_2+2)(e-4)+h-2,   \ h > 0.$
\item  If  $\  r=e-2$, then 
$b=(l_1+l_2+2)(e-1)+h-p-1 $ \ and 
 $\ p \leq 2l_1+2.  $ \\  Further, $p = 2l_1+2 \I h >0.$ 
 \item If $r=e-1$, then 
 $b=(l_1+l_2+2)e+h. $ 

   \end{enumerate} 
 \end{lemma}
 
\underline{Proof}. Recall that  by (\ref{dopoD1}) $e-r\leq 3$  and   by (\ref{set3})    \par $v(R)\!=\!   \{0, e,2e,...,pe, c,\rightarrow\} \cup 
\{y_{1},  ...,y_{1}+l_{1}e\}
\cup 
\{y_{2},  ...,y_{2}+l_{2}e\}$. \\
  Formula  (\ref{XYZ}.1)  with $k=3$ becomes  \par  $b=(r+1)(l_1+l_2+2)-(p+1)(e-r-1)+h$. \\By substituing     $r$ with $ e-1,e-2,e-3 $, we get the desired
expressions for
$b$ in items  3, 2, 1,  respectively. To complete the proof of   (1) and (2), apply  (\ref{l=k}.3) and (\ref{l=k}.2), respectively.\quad  
$\diamond$  
 
\begin{prop} \label{aggiunto1}   Let $ q\in\nat$ be such that  $0<b\leq q(r-1) $. Then:
  \enu
\item    $ r\geq e-q-1$. In particular,  
  \begin{enumerate} 
  \item  If $ r= e-1-q$, then $b=q(r-1), q\leq e-3$,     and
 the equivalent conditions of $(\ref{dopoD}.4)$ hold. 
\item  If $ r\geq e-q$, then  $e-r\leq
k \leq q.$ \end{enumerate} 
\item  If    $(q-1)(r-1)< b\leq q(r-1)$,   we have: 
\enua\item[$(c)$] $k-1 \leq q\leq n-1  $. 
\item[$(d)$] $(q-k-1 )(r-1)< \sum_{i\in
A} (r-r_i )\leq(q-k )(r-1)+e-1-k.$ 
\denu

 \end{enumerate}   
\end{prop}

\underline{Proof}. (1). First we deduce the inequalities   \\ \centerline{ ($\odot$)   \qquad
$(e-r-1)(r-1)\leq b \leq q(r-1),$\qquad} 
by combining (\ref{dopoD}.1) with the assumption. Hence we get   $r \geq e-q-1$.\\  
 ($a)$.  If $ r= e-1-q$, then relations $(\odot)$ give $b=(e-r-1)(r-1)$, and so  the  conditions of  (\ref{dopoD}.4) hold. \\
 ($b)$.  Assertion $(**)$ in the proof of (\ref{dopoD}) insures that $kr-(e-1) \leq b$. Hence   assuming  $r \geq e-q$, we have $ r 
k
\leq b+e-1\leq q(r-1)+e-1\leq  (q+1)r-1$; then $k\leq q$.\\
(2).  Put $M:=\sum_{i\in A}(r-r_i )$. We  have to compare the two inequalities of
(\ref{dopoD}.1)\par $(k-1)(r-1)  +M  +k-(e-r)\leq b\leq k(r-1)+M$ \\  with the assumption\par
$(q-1)(r-1)< b\leq q(r-1)$. \\We obtain the following:  \par
 $(k-1)(r-1)  +M+k-(e-r)\leq q(r-1)$, and also\par
$(q-1)(r-1)  < k(r-1)+M$.\\
The first inequality gives $q\geq k-1$ and \ $M\leq  (q-k )(r-1)+e-1-k$.\\
 The second one says that  $ M> (q-k-1)(r-1).$  Moreover, combining the hypothesis with (\ref{dopoD}.2)  \par $(q-1)(r-1)<b \leq (n-1)(r-1)$, \\ and this implies   $q\leq
n-1 $, as desired.
\quad $\diamond$

 \section {Classification.} Our aim is now to classify    the value sets for one dimensional local domains having
\\
\centerline{$0\leq b\leq 2(r-1)$. }
On this topic several results are present in the literature. For   semigroup rings
$R=\kappa[[t^{\alpha}; \alpha \in S]], \ S\subset \nat$ a numerical semigroup, Brown and Herzog in \cite[Corollary after Theorem 4]{bh}  illustrate the case 
$b =1$.  This result can be  extended to rings $R$ as in Setting \ref{set2} 
(see (\ref{b=0})). Successively    D. Delfino in \cite[Corollary 2.11 and Corollary 2.14]{d} gives a characterization of rings satisfying the condition
$b<r-1$  and    exhibits all the possible value sets in the case $b\leq r $, under the additional
assumption
$r=e-1$. See also Proposition 2.7 from \cite{dd}  for a further generalization. An  exaustive
description of the  cases $0<b< r-1$
 can be found in  \cite[Th.4.6]{ozinv}. \\
 
 In this section we assume Setting
\ref{set2} and Notation
\ref{notaB}. Moreover, 
$t.s.(R)$ will denote the type sequence of
$R$, defined in (\ref{set2}).\\

 First we recall in (\ref{b=0}) and (\ref{finor-1}) the quoted known results,  which now become an easy  consequence of 
our preceding statements. \\

 \begin{thm} \label{b=0}  {\rm Case} $b=0$.\\ The following conditions are equivalent:
\enu  
\item $  b=0 $.
\item  Either $   R$ is Gorenstein, 
  or $ v(R)=\{0,e,..,pe,(p+1)e \rightarrow\}$. 
\item  $t.s.(R)=[r,...,r]$.
\denu
\end{thm}
 \underline{Proof}. By (\ref{dopoD}.1)  $0=b \geq(k-1)(r-1)$; hence either $r=1$ or $k=1$, and this last condition gives, by (\ref{delta1}),
   $ v(R)=\{0,e,..,pe,(p+1)e \rightarrow\}$, or equivalently, $t.s.(R)=[e-1,...,e-1]=[r,..,r]$. Hence $ 1 \I  2 \II  3 $ are clear. Of course, in the Gorenstein 
case
  we have $   t.s.(R)=[1,..,1]$. Implication $ 3 \I  1 $ is immediate by (\ref{n-1}.1). \quad $\diamond$ 
  
\begin{thm}  \label{finor-1} {\rm \cite[Theorem 4.6.1]{ozinv}}   {\rm Case} $0<b<r-1 $.\\ The following facts are equivalent:
 \enu
  \item $ 0<b< r-1. $ 
 \item $v(R)=\{0,e,..,pe,c\rightarrow\}$ \ with \ $pe+2<c \leq (p+1)e.  $  
 \item $t.s.(R)=[ e-1, e-1,...,e-1, r_n], \ r_n>1  $.  
 \end{enumerate}
If these conditions    hold, then:\par$  b<e-2 , \   r=e\!-\!1, \ r_n=e \!-\!1\!-\! b, \ k =1,\
c=(p+1)e-b.  $
\end{thm}

\begin{thm}  \label{b=r-1} {\rm Case} $b=r-1$. \\ If  $b =r-1>0,  $ then   either  $r=e-1$   or $r=e-2$. 
\begin{enumerate}
\item {\rm Subcase}  $r=e-1$. \  The following facts are equivalent:
              \begin{enumerate}
           \item $  \   b  = r-1>0$ \  and  $\   r=e-1 $. 
            \item \ $v(R)=\{0,e,..,pe,pe+2\rightarrow\}, \ e>2$.  
            \item \ $t.s.(R)=[ e-1,....,e-1, 1], \ e>2$. 
           \item  \ $  b  = r-1>0$ \  and  $\  k=1 $. 
            \end{enumerate}

\item {\rm Subcase}  $r=e-2$. \ The following facts are equivalent:
           \begin{enumerate}
             \item [$(e)$] $ b= r-1>0$ \  and  $\ r=e-2 $. 
             \item [$(f)$]  either  \ $v(R)=\{0,e,2e-1,2e,3e-1\rightarrow\}, \ e>3, $  \par
                 or \hspace{0.55cm}   $v(R)=\{0,e, y,2e\rightarrow\}$,   with   $ 2y<  3 e, \ e>3 $.  
           \item  [$(g)$]  either \ $t.s.(R)=[e-2,e-2,1,e-2]$,   with   $ \ e>3,$ \par
                           or \ \ \ \ \ $t.s.(R)=[ e-2,\ r_{2},\ r_3]$,    with $r_{2}\!+r_3\!=e-1, \ e>3$.
 \item[$(h)$]   \  $ b  = r-1>0$ \  and  $\  k=2 $. 
                 \end{enumerate} 
                 
\end{enumerate} 
\end{thm} 

\underline{Proof}. \ Applying  (\ref{aggiunto1}.1) with $q=1$, we obtain that $r\geq e-2$.  Further,  if $b=r-1$, then $r=e-2 \II k=2$ by
(\ref{aggiunto1}.1$a$)  and (\ref{dopoD}.4); also, if    $b=r-1$, then $r=e-1 \II k=1$ by (\ref{aggiunto1}.1$b$). This proves
the first assertion and the equivalences $ a 
\II   d $, $ e   \II  h  $.\par
(1). First note that $(a)$ implies $e>2;$  in fact,  $e=2$   would imply $r=1, b=0$.
$ d 
\I   b $. Since 
$k =1$,   the equivalent conditions of (\ref{delta1}) hold, and\par 
$v(R)=\{0,e,..,pe,c\rightarrow\},$  with $c=(p+1)e-b=pe+2,  e>2$. \\   $ b \I  c $. If $(b)$  holds, then by (\ref{delta1})  $t.s.(R)=[ e-1,....,e-1, r_n]$, with
$r\!-\! r_n=b=r-1$, hence $r_n=1$, as in $(c)$.\\
 $ c \I  a $. By (\ref{n-1}.1),  $   (c)$ implies  $r=e-1$ and $b=r-r_n=r-1$, as in $(a)$.\par
 (2). First note that $(a)$ implies $e>3;$  in fact,  $e=3$   would imply $r=1, b=0$.\\ $ h \I  f $. Since  
 $ k= 2  $   we use (\ref{Casek=2}.2) recalling that   $ p \leq 2l+2$:\par
 $e-3=b=(l+1)(e-1)+h-p-1 \geq (l+1)(e-1)+h-2l-3$.\\
Hence we get  $l(e-3)+h \leq 1$  and the following possibilities occur by (\ref{Casek=2}.2$c$):\par
 $ (l,p,h)=(0,1,0),$ \  or  \ $  (l,p,h)=(0,2,1),$  \ or \ \ $h=0,\  e=4,\   l=1$. \\
 - If $  (l,p,h)=(0,1,0)$, then $   c=2e$,     $v(R)=\{0,e, y,2e\rightarrow\}$  with  $ 2y<  3 e, \ e>3$.\\
 - If $ (l,p,h)=(0,2,1)$, then   $ v(R)=\{0,e,2e,c \to \} \cup \{y\}$,  with
$c-\delta=4$,\par  $ c=(p+1)e-h=3e-1, \ 2y<c+e=4e-1 \I y \leq2e-1,$\par $  c-e \in v(R)  \I y=2e-1$. Hence 
$v(R)=\{0,e,2e-1,2e,3e-1\rightarrow\}, \ e>3 $. \\
-  If $  h=0, \ e=4,  \ l=1, $ then $ e-3=b=(l+1)(e-1)+h-p-1 \ \I$\par$  p=4  =2l+2 \I h>0$, which is absurd. 
Hence $ h   \I  f $ is proved. \\
  $ f  \I  g $. Denoting  $R_0=\kappa[[t^d, \ d\in v(R)]]$   the monomial ring such that  
$v(R_0)=v(R)=  \{0,e,2e-1,2e,3e-1\rightarrow\} $, we have $r(R_0) =e-2$.  Since $r(R) \leq
r(R_0)$   and
   $r(R) \geq e-2$   by   (\ref{dopoD1}.1), we conclude that    $r(R)  =e-2$.   The other  invariants are easily derived from     $v(R)$: 
$ c-\delta=4,\ \delta=3e-5,\ b=(c-\delta)r-\delta=e-3$.   By substituting in    
(\ref{dopoD}.1),  we obtain   $\sum_{h\in A} (r-r_h)=0$,  hence   $r_2=e-2$  
and  $\ r_3+r_4=e-1$, as desired.  \\
The same reasoning holds for $v(R)=\{0,e, y,2e\rightarrow\}$. \\
To see   $ g  \I  e $, it suffices to recall that   $b= \sum_{h=1}^n(r-r_h)$, see  (\ref{n-1}.1). \quad $\diamond$

\begin{thm} \label{b-compreso}  {\rm Case} $  r-1 < b < 2(r-1)$.  \\ We have \ $  r-1 < b < 2(r-1)$ if and only if $ v(R)$ is
one of the following:
\enu
\item $ v(R)=\{0,e,....,pe,c\rightarrow \}\cup \{y\},$  with  \   $ y \notin e \Bbb Z$,  \\    and either   $ 2y\geq c+e, \   pe+5\leq c \leq min\{y+e,(p+1)e\},  \  e\geq 5, $\\
or  \ $e=2e', \  y=3 e', \  p=2, \ 4e'+5\leq c \leq 5e', \ e \geq 10, \ y\in v( xR:\m).$\ 
 
\item $ v(R)=\{0,e,2e,c\rightarrow \}\cup \{y\},$   with   $\ y \notin e \Bbb Z, \ 2y< c+e$ and:\\ if $2y \neq 3e$, then  $\ 2e+3\leq c \leq  3e-2, \ e\geq 5;$ \\ if   
$2y=3e$,   then $ e=2e',\ 4e'+3 \leq c \leq 5e', \ e \geq 6, \ y\notin v( xR:\m).$\ 
  
\item $  v(R)=\{0,e,y,c\rightarrow \},$ \\  with   $ y \notin e \Bbb Z, \ e\geq 5, \  2y<   c+e, \ e+4\leq  c \leq 2e-1.$  
 \end{enumerate} 
 In each case $k=2$; in case $(1)$,  $r=e-1$ and   $b\geq r+1$; in cases $(2)$ and $(3)$,   $r=e-2$.  
\end{thm}

\underline{Proof}.  Assume  $ r-1<b<2(r-1)$. \par
Step 1.  Claim:  $k=2$ \   and  \  $e-2\leq r\leq e-1, \ r>2$.\\ We have  $ r>2$, since $r=2\I 1<b<2,$ which is absurd. 
Further 
(\ref{dopoD}.1) gives  $(k-1)(r-1) \leq b$, and so
$k\leq 2$. But   $k=1$ would imply $b\leq r-1$ by (\ref{delta1}), then
$k=2$. We conclude using (\ref{dopoD1}.1).
\\ Now utilizing the notation in  (\ref{set3}) we   write:\\
$
(*)\left\{\begin{array}{lll}
\!\! v(R)\!=\!   \{0, e,2e,...,pe, c,\rightarrow\} \cup  
\{y ,y +e,...,y +l e\}, \  p \geq 1, \ y>e, \  y \notin e \Bbb Z,\\
 \!\! y +l e <c=(p+1)e-h \leq y +(l +1) e,\ \ l +1\leq p.
\end{array}\right. $

   Step 2.  Claim: $\ l=0$ \  and  \ $e\geq 5$. Further, if $r=e-2$, then $\ p\leq 2$.\\ 
- If $r=e-1,$ then,  by (\ref{Casek=2}.1) we know that $b=(l+1)e+h $, \ $l, h\geq 0$. Hence $b<2(r-1) =2e-4 \I (l-1)e+h <-4 \I$  $ l=0,\  h<e-4;$ \
further we get : $c=(p+1)e-h
\geq pe+5 $,\ \  $e\geq 5$ and $b= h+e \geq e=r+1$.\\  
- If $  r=e-2,$  we have     $\ (l+1)(e-3)
\leq  b $ \ and  $\ l+1\leq p\leq 2l+2 $ \  by (\ref{Casek=2}.2). Then
$b<2(r-1)=2(e-3)\I l=0$ and $p\leq 2;$  also, the assumption  $e-3<b<2e-6$ implies $ e\geq 5$. 
\par Step 3. When $ r=e-1,$  recalling the relations proved
in Step 2,  we obtain  
$ v(R)=\{0,e,....,pe,c\rightarrow \}\cup \{y\}$, with $ e\geq 5,  \  pe+5\leq c$, as in item 1. Recall that by definition of $p$ and $l$ we have $c \leq (p+1)e$ and $c \leq y+e$. Moreover, by
 (\ref{l=2}.1) one of the following conditions is satisfied: \\ either  $(a)$
  $ 2y\geq c+e,$  or $(b)$  $2y=(2q+1) e<c+e,\ p\geq 2$  and $y\in v( xR:\m).$ \\ Further,  as noted in  (\ref{l=k}.4), $p\geq 3, \ l=0 \I 2y >c+e$, hence in
case 
$(b)$ we have $p=2$ and consequently $(2q+1) e<c+e\leq 4e\I q=1$. This proves (1).
\par Step 4. When $ r=e-2,$ we have by Step 2 that $ l=0$ and $p\leq 2$. \\ In the case $p=2$ we get item $2$. In fact from (\ref{Casek=2}.2)  we
obtain  $c= 4e-4-b $    and the bounds for $c$ follow at once. The last assertion in item
$ 2 $ comes from (\ref{l=2}).  Analogously, in the case $p=1$ we get item $3$. Notice
that when
$p=1$ we cannot have
$2y=3e< c+e$, because
$c+e\leq 3e-1$.  
\par

To complete the proof, let $v(R)$ be as in items $(1), (2), (3);$   we claim that  $ r-1 < b < 2(r-1)$. In every case $k=2$; in case (1) $r=e-1$ and in cases (2), (3)  
$r=e-2$ by (\ref{l=2}). The rest is a direct computation based on relation (\ref{XYZ}.2): $c=(p+2)(r+1)-b$. 
\quad $\diamond$

\begin{example}{\rm  We supply an example for each case of the above proposition. \par$\bullet$ Case (1) with $2y\geq c+e$. \\Let $R=\kappa[[t^{5}, t^{10}, t^{12}, t^{15},\to]]$. Then:
$y=12,   p=2,   c=15,   r=4,  b=5 $.\par
 $\bullet$ Case (1) with $2y=3e$.\\
 Let $R=\kappa[[t^{10}, t^{15}, t^{20}, t^{25},\to]]$. Then: $y=15, p=2, c=25, r=9, b=15$.\par
 $\bullet$ Case (2).\\ 
  Let $R=\kappa[[t^{10}, t^{15}+t^{16}, t^{20}, t^{25},\to]]$. As above, $y=15, p=2, c=25$, but $r=8$ by \ref{tipoe-1} since 
 $(t^{15}+t^{16})^2 \notin x \m$. Then $b=11$.
 \par $\bullet$ Case (3).\\ Let $R=\kappa[[t^{5},t^{6}, t^9 ,\to]]$. Here $y=6, p=1, c=9, r=3, b=3  $.   } \end{example}

  \begin{thm} \label{b=2(r-1)}  {\rm Case $ b= 2(r-1)$.}  \\ $ b=2(r-1)>0$ if and only if $ v(R)$ is one of the following:
\begin{enumerate}
\item  
 \begin{enumerate}
 \item  $v(R)\!=\! \{0, e, e+2, e+4\rightarrow\}$, \      $e\geq 4$.
 \item $v(R)\!=\! \{0, e,  2e, 2e+4\rightarrow\}\cup\{y\},$ \     $e\geq 4, \ y \in v(xR: \m)$.
\item $v(R)\!=\! \{0, e, 2e,..., pe, pe+4\rightarrow\}\cup\{y\}$, \   $e\geq 4, \ y \geq (p-1)e+4, \ p\geq 3$.
 \end{enumerate}
\item    
 \begin{enumerate}
   \item $v(R)\!=\!\{0,e,e+1,  e+3\rightarrow \} $, \  $e\geq 4$. 
 \item $v(R)\!=\! \{0, e, y, 2e, 2e+2\rightarrow\},$ \    $e\geq 5$, \ $2e+4\leq 2y <3e+2, \ 2y \neq 3e$. 
 \item $v(R)\!=\!\{0,e, 2e,3e-1,3e, 4e-1, 4e,  5e-1 \rightarrow \}$, \  $e\geq 4$.
\item $v(R)\!=\!\{0,e, 2e,y, 3e, y+e, 4e \rightarrow \}$,  \   $e\geq 4$,  \ $2y<5e $.
 \end{enumerate}
\item  
 \begin{enumerate}
 \item $v(R)\!=\! \{0, e, y_1, y_2,2e \rightarrow\}$, \  $e\geq 5$, \  $y_1+y_2<3e $.
\item $v(R)\!=\! \{0, e,2e-2,2e-1,2e,3e-2 \rightarrow\},  $ \  $e\geq 5$.
 \end{enumerate}
 \end{enumerate} 
Further:\par
in case  1, $\quad r=e-1$ and $ \ell_R (R/({\goth C}+xR))=2;$\par
  in case 2, \quad$ r=e-2 $ and $ \ell_R (R/({\goth C}+xR))=2;$\par
    in case 3, \quad$  r=e-3 $ and $  \ell_R (R/({\goth C}+xR))=3$.
\end{thm}

\underline{Proof}. Let, as above, $k=  \ell_R (R/({\goth C}+xR))$. First we assume  $b= 2(r-1)>0$ and we observe  that by (\ref{dopoD}.1) $(k-1)(r-1) \leq b=2(r-1)$, then $k \leq 3$. Since $k=1$ implies $b \leq r-1$ by   (\ref{delta1}), one of the following cases occurs:\par

$\left [\begin{array}{l}or\quad  k=2 \quad and \quad r=e-1  \\
or\quad  k=2 \quad and \quad  r=e-2 \\
or\quad k=3 \quad and \quad r=e-3.\end{array}\right.$ \\
In case  $k=2$ by Setting  \ref{set3}    we have:\par
$(*)\left\{\begin{array}{lll}
\!\! v(R)\!=\!   \{0, e,2e,...,pe, c,\rightarrow\} \cup  
\{y ,y +e,...,y +l e\}, \  p \geq 1, \ y>e, \  y \notin e \Bbb Z,\\
 \!\! y +l e <c=(p+1)e-h \leq y +(l +1) e,\ \ l +1\leq p.
\end{array}\right. $

Step 1. Assuming  $  r=e-1$ and $k=2$,   we   prove that $v(R)$ has the form described in item 1. 
By (\ref{Casek=2}.1) and the assumption   we have the equalities
 $ b =  (l+1)e+h=2e-4  $;   hence 
 $ (l-1)e+h=-4 \I l=0, \ h=e-4, \ e \geq 4, \ c= (p+1)e-h=pe+4$.  
Now, $l=0\I   y \geq c-e=(p-1)e+4$, and so  \par
 $ v(R)\!=\! \{0, e,..., pe, pe+4\rightarrow\}\cup\{y \}$,  \  with $\ (p-1)e+4\leq y  \leq  pe+2,\ \ e\geq 4. $ 
For $p=1$ we get   (1.$a$). In fact, by (\ref{l=2}.1)  $2y \geq c+e=2e+4 \I y \geq e+2\I y=e+2$. 
For $p=2$ we get  (1.$b$).
For $p \geq 3$ we get (1.$c$).\par

Step 2. Assuming   $  r=e-2$  and $k=2$,   we   prove that $v(R)$ satisfies item 2. 
First, 
by (\ref{Casek=2}.2) we have that $l+1\leq p\leq 2l+2$ and also that\par  $(**)  \quad (l+1)(e-3)\leq (l+1)(e-1)+h-p-1=b $. \\Then $b=2(e-3) >0$ implies $(l+1)(e-3) \leq
2(e-3)$, i.e.  
$l\leq 1.$
\\ Case $l=0$, and consequently $ 1\leq p\leq 2$.\par
$(\cdot )$  If $l=0, \ p=1,$  then by $(**      )$,  $h=e-3$, thus $c=e+3$, and $ (2.a) $ holds.

$(\cdot )$   If $l=0, \ p=2,$  then $h=e-2,   c=2e+2$, hence $ (2.b) $ holds.\\
Case $l=1$. Now, relation $(**)$ combined with the assumption $b=2e-6$ implies \\ $ h-p-1=-4,\  2\leq p\leq 4 $ and   two possibilities occur:\par

 $(\cdot )$  $p=4,  \ h=1, \ c=5e-1$.  The relation $c \leq y +(l+1)e$ gives    $y  \geq 3e-1,$ the relation $2y <c+e$ gives
$y  \leq 3e-1.$   Hence    $ (2.c)  $ holds.

 $(\cdot )$ $p=3, \ h=0, \ c=4e$; hence $ (2.d) $ holds. \par
 
 Step 3. Assuming    $  r=e-3$   and $k=3$,   we   prove that $v(R)$ has the form described in item 3. 
 First, by Setting \ref{set3} and by (\ref{Casek=3}.1)  we have:\\
$ (\overline *)\left\{\begin{array}{lll}
\!\! v(R)\!=\!   \{0, e,...,pe, c\rightarrow\} \cup  
\{y_1 ,y_1 +e,...,y_1 +l_1 e\}\cup\{y_2 ,y_2 +e,...,y_2 +l_2 e\} \\ 
 p \geq 1, \ y_2>y_1>e, \  y_i \notin e \Bbb Z, \\
 \!\! y_i +l_i e <c=(p+1)e-h \leq y_i +(l_i +1) e,\ \ l_i +1\leq p,\\
b=(l_1+l_2+2)(e-2)+h-2(p+1).
\end{array}\right. $
 By (\ref{l=k}.3), since $r=e-k$, then \ $p \leq l_1+l_2+2$. 
 
  $(\cdot)$ If  $p<l_1+l_2+2$,  then  substituting     $b=2(e-4)>0$  in  $ (\overline *)$  we get  $   (l_1+l_2)(e-4)+h \leq 0,   \ h \geq 0.$   Hence    $h=l_1=l_2=0, \
  p=1, \ c=2e, \ y_1+y_2<c+e $ by (\ref{l=k}.3), and so we have  $(3.a)$.
  
$(\cdot)$   If  $p=l_1+l_2+2$,  then  analogously we get $ (l_1+l_2)(e-4)+h-2=0,$ with $0<h   \leq 2$. The case $h=1$ is impossible. In fact,  $h=1 \I l_1+l_2=1$ (in particular, by (\ref{aftersetting}.2), $l_2 \leq l_1$, hence $l_2=0, \ l_1=1 ), \  e=5, \
  p=3, \
 c=(p+1)e-h=19$. The relation of (\ref{set3}) $c \leq y_i+(l_i+1)e$ gives $y_1\geq 19-10=9, \ y_2 \geq 19-5=14$, but $y_1+y_2 <c+e=24$ by (\ref{l=k}.3); the only possibility would be 
$ y_1=9, \ y_2=14.$ Absurd that $\overline{y_1}=\overline{y_2} \ (mod \ 5)$.  Hence    $h=2, \ l_1=l_2=0, \ p=2, \ c=3e-2 $ and \par
 $ v(R)\!=\! \{0, e,2e, 3e-2,\rightarrow\}\cup \{y_1,y_2\}$. \\
Since $l_1=0$, the bound   $c \leq y_1+e$  gives  $y_1\geq 2e-2$.  Recalling that  by (\ref{l=k}.3)
 $y_1+y_2< c+e$,  we conclude $ y_1=2e-2,   y_2=2e-1$, as in $(3.b)$.  \par
 
  Viceversa, we assume in the following $v(R)$ having the form described in items 1,2,3, and we prove that   $b= 2(r-1)>0$.\\
  For a $v(R)$ as in item 1 we see that $r=e-1$ using (\ref{l=2}). In fact, in case $(1.a)$ we have $y=e+2, \ 2y=c+e$ and in case $(1.c)$:\par  $2y \geq 2(p-1)e+8>c+e=(p+1)e+4$.\\
  In conclusion in each case of item 1 we have  $ \ell_R (R/({\goth C}+xR))=2, \ r=e-1, \ l=0$. Using (\ref{Casek=2}.1) $b=e+h=2e-4=2(r-1)$, as desired.\\
  In case $(2.a)$,  $y=e+1 \notin  v(xR: \m)$, then $r=e-2$ by (\ref{tipoe-1}).
  In case  $(2.b)$ by hypothesis $2y<c+e$ and $2y \neq 3e$, then  $r=e-2$ by (\ref{l=2}).
 In case  $(2.c)$ we get by a direct calculation $ v(xR_0: \m)\setminus v(R_0)=\{4e+1,....,5e-2\}$, then $r=r(R_0)=e-2$. 
 In case  $(2.d)$ $2y<c+e$ and $2y \notin e \Bbb Z$, then $r=e-2$ by  (\ref{l=2}). 
  In conclusion, in each case of item 2 one has: $ \ell_R (R/({\goth C}+xR))=2, \ r=e-2$, and so by (\ref{Casek=2}.2) $b=(l+1)e-1)+h-p-1$. Putting in this formula\par
$(\cdot )$ $l=0, \ p=1, \ h=e-3,$ in case $(2.a)$,

$(\cdot )$ $l=0, \ p=2, \ h=e-2,$ in case  $(2.b)$,

$(\cdot )$ $l=1, \ p=4, \ h=1,$ in case  $(2.c)$,

$(\cdot )$ $l=1,
       p=3, \ h=0,$ in case  $(2.d)$,\\ we get $b=2e-6=2(r-1)$, as desired.\\
  In both cases of item 3 we have $r=e-3$. In fact, $y_1+y_2-e \notin v(R) \I y_1, y_2 \notin v(xR: \m) \I e-r=  3$ by (\ref{l=k}.1).   
 Hence   $ \ell_R (R/({\goth C}+xR))=3, \ r=e-3, \ l_1=l_2=0$, and   by (\ref{Casek=3}.3)  $b= 2 (e-2)+h-2(p+1)$.   Putting in this formula\par
    $(\cdot )$  $h=0, \ p=1$ in case  $(3.a)$,
   
     $(\cdot )$ $h=p=2$ in case  $(3.b)$, \\
      we get $b=2e-8=2(r-1) $, as desired. $\quad \diamond$\\
      
 \noindent With similar arguments one can evaluate the semigroups $v(R) $ of rings having $b>2(r-1)$. For instance, if  $2 (r-1)<b\leq 3(r-1)$  there are few possible cases and the
classification is tedious but easy. Now, for each $q\geq 3$ we  construct a   family of rings   of multiplicity $e$ and Cohen Macaulay type $r=e-1$ having $b=q(r-1)$
 or $(q-1)(r-1)<b<q(r-1)$.  
 
\begin{example}   {\rm   Let  $q\geq 3$.  Following notations of Setting \ref{set3} we consider \par
 $v(R)=\{0,e,2e,...,pe, c \to \} \cup \{y,y+e,...,y+le\}$, \\
with $e >p, \  p=2q, \ l=q-2$.  In this case $k=2$. Using (\ref{l=2}) we see that
$r=e-1$, because $y+(q-1)e \geq c >2qe \I y>(q+1)e \I 2y >2(q+1)e\geq c+e$. Then by (\ref{Casek=2}.1) $b=(q-1)e+h$, with $0 \leq h \leq e-2$.   Now, with an  additional hypothesis on the
conductor, we are in goal. In fact:  

 1) Assuming $ c=pe+p$, we have $h=(p+1)e-c= -p+e=-2q+e$, then $b=(q-1)e+( -2q +e)=q(e-2)=q(r-1)$.  
 
 2)  Assuming $ c> pe+p$, i.e.  $e-h > 2q$, we have
  $(q-1)(e-2) < (q-1)e \leq b=(q-1)e+h= q(e-2)+2q-e+h < q(e-2)$, hence   $(q-1)(r-1)<b<q(r-1)$.}\end{example}

As a further application of the previous results we describe exhaustively the cases $b=1$ and $b=2$ (see next (\ref{b=1}), (\ref{b=2}); for $b=1$ see also \cite{bh}, Section 4). With regard  to the formula
$$ b=\sum_{i=1}^n (r-r_i) $$
it becomes natural to consider the invariant $b$as a   measure of how far is the type sequence $[r_1,...,r_n]$ from the maximal one $[r,...,r]$. For instance, for $b=1$ one expects a type sequence of the form $[r,...,r-1,...,r]$, for $b=2$ $[r,...,r-1,...,r-1,...,r]$ or $[r,...,r-2,...,r]$, and so on. Surprisingly, after finding by a direct computation all the possible value sets and the corresponding type sequences, we discover that very few choices are possible. For  $b=1$ (resp. $b=2$) either $e \leq 4$ (resp. $e\leq 5$) or $t.s.(R)= [e-1,...,e-1,e-1-b]$.

 \begin{coro} \label{b=1}  \ {\rm Case} $b=1$. \
  Here  t.s. stands for $t.s.(R)$. \\$ b=1 $ if and only if $ v(R)$ is one of the following:\par
 $v(R)=\{0,4,7,8,11\to\}$, with
t.s.   $[2, 2,1, 2];$ \par $v(R)=\{0,4,5,8, \to\}$,  with t.s. $[2,1,2];$  \par    $v(R)=\{0,e,...,pe,(p+1)e-1,\to\}, \ e \geq 3, $ with t.s. $[e-1,...,e-1,e-2]$.
 \end{coro}

\underline{Proof}. First recall that $b>0 \I r>1$ by (\ref{n-1}.1). Let, as in  (\ref{XYZ}.1), 
$ b=  X+ Y +Z$, 
 where  \ $X:=(k-1)(r-1)\geq 0  , \ Y:=k-(e-r) \geq 0$, \ and \\ $Z:=(r+1)(p+\sum_1^{k-1}  l_i)+k+h-pe-1\geq 0$.\\ 
Assuming $b=1$,   we have to consider the choices: \\
\centerline{ $ \matrix{\ &X &Y&Z \cr
a)&1&0&0\cr
b)&0&1&0\cr
c)&0&0&1}$ }
 In $a)   \ k=r=2, \ 2-(e-2)=0 \I e=4$. By  (\ref{b=r-1}.2) with $e=4$ we find: \par   $v(R)=\{0,4,7,8,11\to\}$,   \par  $v(R)=\{0,4,5,8, \to\}$. \\
 In $b) \ k=1, \ 1-(e-r)=1$,  which is absurd.\\
 In $c) \ k=1, \ 1-(e-r)=0 \I r=e-1, \ e \geq 3, \ Z= ep+1+h-pe-1=1 \I h=1 \I c=(p+1)e-1$. By  (\ref{delta1})   we find: \par  
$v(R)=\{0,e,...,pe,(p+1)e-1,\to\}, \ e \geq 3 $.
  $\quad \diamond$

 \begin{coro} \label{b=2}  \  {\rm Case} $b=2$. \
 As above, t.s. stands for $t.s.(R)$. \\$ b=2 $ if and only if $ v(R)$ is one of the following:\par
 $v(R)=\{0,4,5,7,\to\}$,  with t.s.  $[2, 1, 1];$ \par $v(R)=\{0,4,8,11,12,15,16,19,\to\}$,  with t.s. $[2,2,2,1,2,1,2];$\par   $v(R)=\{0,4,8,9,12,13,16, \to\}$, with t.s. $[2,2,1,2,1,2];$ \par
$v(R)=\{0,5,9,10,14,\to\}$,  with t.s. $[3,3,1,3];$ \par 
 $v(R)=\{0,5,6,10,\to\}$, with t.s. $[3,1,3];$\par
   $v(R)=\{0,5,7,10,\to\}$,   with t.s. $[3,2,2];$ \par 
 $v(R)=\{0,5,6,7,10,\to\}$,   with t.s. $[2,1,1,2];$\par 
  $v(R)=\{0,5,6,8,10, \to\}$,  with t.s. $[2,2,1,1];$ \par 
 $v(R)=\{0,5,8,9,10,13, \to\}$,   with t.s. $[2,2,1,1,2];$  \par
  $v(R)=\{0,e,...,pe,(p+1)e-2,\to\}, \ e \geq 4$,  with t.s. $ [e-1,...,e-1,e-3] $. 
 \end{coro}

\underline{Proof}. 
 As in the preceding proof, assuming $b=2$,   we have to consider the following choices: \\
 \centerline{$\matrix{\ &X &Y&Z \cr
a)&0&1&1\cr
b)&1&0&1\cr
c)&1&1&0\cr
d)&2&0&0\cr
e)&0&2&0\cr
f)&0&0&2}$}
 First recall that $k=1 \I r=e-1$ by (\ref{delta1}), and so $X=0$ (with $r>0$) $\I k-(e-r) =Y=0$ and cases $a), e)$ are impossible.\\
  In $b)  \ X=1 \I k=r=2, \ 2-(e-2)=Y=0 \I e=4 $, hence $b=2(r-1)$ and we can apply   (\ref{b=2(r-1)}.2) with $e=4$. We find: \par  $v(R)=\{0,4,5,7,\to\}$,  \par  $v(R)=\{0,4,8,11,12,15,16,19,\to\}$, \par  $v(R)=\{0,4,8,9,12,13,16, \to\}$.\\ 
 In $c) \ X=1 \I k=r=2, \ 2-(e-2)=Y=1 \I e=3, \ Z=3(p+l)+2+h-3p-1=0 \I 3l+h+1=0$, which is absurd.\\
 In $d)$ the condition $X=(k-1)(r-1)=2$ implies two possibilities: \\
 $d_1)   \ k=2, \ r=3, \ 2-(e-3)=0 \I e=5$.  We are in case $b=r-1, \ r=e-2$. By   (\ref{b=r-1}.2) with $e=5$ we find: \par 
 $v(R)=\{0,5,9,10,14,\to\}$,\par 
  $v(R)=\{0,5,6,10,\to\}$, \par
 $v(R)=\{0,5,7,10,\to\}$.\\ 
 $d_2) \ k=3, \ r=2, \ e=5 $.  We are in case $b=2(r-1), \ r=e-3$, and so by   (\ref{b=2(r-1)}.3) with $e=5$ we find: \par 
 $v(R)=\{0,5,6,7,10,\to\}$,\par 
 $v(R)=\{0,5,6,8,10, \to\}$,\par 
  $v(R)=\{0,5,8,9,10,13, \to\}$.\\ 
 In $f) \ k=1, \ r=e-1, \ Z=ep+1+h-pe-1=2 \I h=2 \I c=(p+1)e-2$. By  (\ref{delta1})   we find: \par  $v(R)=\{0,e,...,pe,(p+1)e-2,\to\}, \ e \geq 4 $.  $\quad \diamond$

\end{document}